\documentclass[times,review]{elsarticle}

\usepackage{subfigure}
\usepackage{jcomp}
\usepackage{framed,multirow}

\usepackage{amssymb}
\usepackage{latexsym}

\usepackage{url}
\usepackage{xcolor}
\definecolor{newcolor}{rgb}{.8,.349,.1}

\journal{Journal of Computational Physics}

\begin{document}

\verso{W. Chen \textit{et al.}}

\begin{frontmatter}

\title{An explicit and non-iterative moving-least-squares immersed-boundary method with low boundary velocity error}%

\author[1]{Wenyuan Chen}
\author[1]{Shufan Zou}
\author[1]{Qingdong Cai}
\author[1,2]{Yantao Yang \corref{cor1}}
\cortext[cor1]{Corresponding author.}
\ead{yantao.yang@pku.edu.cn}

\address[1]{State Key Laboratory for Turbulence and Complex Systems, Department of Mechanics and Engineering Science, College of Engineering, Peking University, Beijing 100871, China}
\address[2]{Beijing Innovation Center for Engineering Science and Advanced Technology, Peking University, Beijing 100871, China}

\received{}
\finalform{}
\accepted{}
\availableonline{}
\communicated{}

\begin{abstract}
In this work, based on the moving-least-squares immersed boundary method, we proposed a new technique to improve the calculation of the volume force representing the body boundary. For boundary with simple geometry, we theoretically analyse the error between the desired volume force at boundary and the actual force given by the original method. The ratio between the two forces is very close to a constant. Numerical experiments reveal that for complex geometry, this ratio exhibits very narrow distribution around certain value. A spatially uniform coefficient is then introduced to correct the force and fixed by the least-square method over all boundary markers. Such method is explicit and non-iterative, and can be easily implemented into the existing scheme. Several test cases have been simulated with stationary and moving boundaries. Our new method can reduce the residual boundary velocity to the level comparable to that given by the iterative method, but requires much less computing time. Moreover, the new method can be readily combined with the iterative method and further reduces the residual boundary velocity.
\end{abstract}


\end{frontmatter}


\section{Introduction}

Since first proposed by Peskin~\cite{peskin1972}, the immersed-boundary method (IBM) has shown huge advantages in simulating fluid-structure interaction (FSI) and multiphase flows mainly due to its simplicity in dealing with complex geometry and moving boundary. The key idea of IB method is modelling the boundary by a virtual volume force, which is so determined that the resulting flow field satisfies the corresponding boundary conditions. By doing so, the Eulerian meshes over which the governing equations are discretized do not need to conform to the boundary. The efficiency and accuracy of the specific IB method, therefore, heavily rely on how the virtual volume force is determined, and on how the calculated force is applied to the Eulerian meshes. Numerous methods have been developed regarding these two aspects. e.g. see the reviews of~\cite{mittal2005,huang2019}.

For a large category of IB methods, namely the diffused interface method, the virtual force representing the boundary is distributed over the neighbouring Eulerian grid points according to certain distribution functions. For instance, a regularized delta function was adopted for the force spreading by Peskin~\cite{peskin1972,peskin2002}. Furthermore, if the boundary condition of the immersed surface is satisfied by imposing kinematic constraints on the control markers on the surface, the method falls into the direct Lagrangian forcing type~\cite{uhlmann2005,vanella2020}. The direct Lagrangian forcing IB method has been applied to various flows, such as hemodynamics~\cite{spandan2017,mittal2016}, multiphase flows~\cite{kempe2012,uhlmann2014}, and turbulent flows~\cite{iaccarino2003}. 

Beside the explicit method, other strategies are employed to further improve the accuracy of the boundary condition, such as the implicit method~\cite{wu2009,wang2011,su2007} and iterative method~\cite{kempe2012,wang2008,ji2012}. The iterative IB method is especially beneficial for flows where the dynamics of the interface requires high accuracy for the velocity field near the immersed boundary~\cite{lau2018,shin2012,calderer2014}. Also for internal flows or flows with closed boundary, to achieve small volume leakage through the boundary, either a very fine Lagrangian mesh must be used, or special treatment should be implemented~\cite{huang2012,le2009,shoele2010}.

In direct Lagrangian forcing method, since the volume force is calculated on the Lagrangian markers on the immersed boundary, the interpolation from Eulerian grids to the Lagrangian markers and the spreading from the Lagrangian markers to the neighbouring Eulerian grids are inevitably involved. Usually some kernel functions are utilized to construct such schemes, such as the regularized delta function~\cite{peskin2002}. The specific design of this delta function has a profound effect on the numerical performance~\cite{shin2008,yang2009}. Various kernel functions have been proposed in past, including the moving-least-squares (MLS) method~\cite{vanella2009}, the reproducing kernel particle method (RKPM)~\cite{pinelli2010}, and inverse distance interpolation~\cite{krishnan2017}, to name a few.

It has been noticed by several groups that, in the direct Lagrangian forcing method, after the spreading and interpolation operations based on certain transfer functions, the actual volume forces added to the Lagrangian markers do not equal to the desired values needed for the boundary condition to be properly imposed~\cite{kempe2012,wang2011}. This is also the reason why implicit or iterative technique was developed. In the present study, we propose an explicit and non-iterative technique for the direct Lagrangian forcing IB method, which is based on the MLS-IBM. The MLS approach itself is one type of meshless approximation~\cite{lancaster1981} and has been implemented into IB method~\cite{vanella2009,tullio2016}. The MLS-IBM is of great advantages for moving and deforming boundaries~\cite{vanella2020,spandan2017,tullio2016}. Usually, iteration is needed to reach a satisfactory level of error on the immersed boundaries, especially for internal flows. Here we aim at achieving similar or higher level or accuracy for boundary condition while avoiding iteration, and therefore saving computing time. 

The rest of the paper is organized as follows. In Section 2 describes the numerical methodology including a detailed description of flow solver and the MLS-IBM. Then we analyse the causes of errors and propose a new numerical method to reduce errors in Section 3. In Section 4 we present computed results for a variety of cases that are intended to validate the solver and to firmly establish its accuracy. Finally, conclusions are presented in Section 5.
 
\section{Governing equations and the baseline method}\label{sec:solver}

\subsection{The fluid phase}

Consider the Navier-Stokes equations and the continuity equation for the incompressible flows
\begin{eqnarray}
   && \partial_t \mathbf{u} + {\mathbf{u}}\cdot{\nabla \mathbf{u}}
     = -\nabla p + \nu\nabla^2\mathbf{u} + \mathbf{f}_b, \label{eq:ns} \\
   && \nabla \cdot \mathbf{u}=0. \label{eq:continua}
\end{eqnarray}
Here, $\mathbf{u}$ is velocity vector, $p$ is pressure, and $\nu$ is kinematic viscosity, respectively. $\partial_t$ denotes the partial derivative with respect to time. $\mathbf{f}_b$ is the volume force introduced in IB method to enforce the boundary condition at the immersed surfaces. 

The governing equations~(\ref{eq:ns}) and (\ref{eq:continua}) will be numerically solved in a Cartesian box. The numerical solver employs a second-order conservative finite-difference scheme with the fractional time-step method on a staggered Eulerian grid. The non-linear terms are discretized by the explicit Adams-Bashforth scheme, and the viscous terms are by the implicit Crank-Nicholson scheme, respectively. For the time integration we use a self-starting Runge-Kutta (RK) scheme with third-order accuracy. An intermediate velocity $\tilde{\mathbf{u}}$ is first computed by considering the non-linear term and the viscous term. Then the divergence-free condition is ensured by solving a Poisson equation for a pressure correction of $\tilde{\mathbf{u}}$. Details of the flow solver found in Refs~\cite{rai1991,verzicco1996,van2015}, along with various validation simulations.

\subsection{The baseline moving-least-squares immersed-boundary method}\label{sec:baseline}

For the volume force $\mathbf{f}_b$ in Eq.~(\ref{eq:ns}), we adopt the MLS-IBM as the baseline method. MLS-IBM has been used by many groups, e.g. see Refs~\cite{vanella2009,kempe2012,tullio2016,spandan2017}. For a detailed discussion of the method the reader is referred to the book chapter by Vanella and Balaras~\cite{vanella2020}. Here we only describe the key ingredients relevant to the current study. In MLS-IBM, the immersed boundary is discretized into triangles, each of which corresponds to a Lagrangian marker with coordinate $\mathbf{X}$. Hereafter, letters with upper case denote the quantities at the Lagrangian markers, and those with lower case denote the values at the Eulerian grids, respectively. 

To impose the boundary condition at the immersed surfaces, a volume force $\mathbf{F}(\mathbf{X})$ is calculated at the Lagrangian markers, which requires the interpolation from a stencil of Eulerian grids $\mathbf{x}_k$ with $k=1,...,ne$ to the Lagrangian point $\mathbf{X}$. Here $ne$ is the number of the Eulerian grids within the stencil. The force $\mathbf{F}(\mathbf{X})$ is distributed to the neighbouring Eulerian points, then the volume force $\mathbf{f}$ at $\mathbf{x}_k$ is the summation over all the evolved Lagrangian markers $\mathbf{X}^l$ with $l=1,...,nl$. $nl$ is the total number of Lagrangian markers which contribute to the volume force at Eulerian grid $\mathbf{x}_k$. The mathematical formula for this procedure read
\begin{eqnarray}
  \mathbf{U}^{L}\left(\mathbf{X}^{l}\right) 
    &=& \sum_{k=1}^{ne} \, \tilde{\bf{u}}_k \, 
        \phi_{k}^{l}(\mathbf{x}_k,\mathbf{X}^l), \label{eq:interp} \\
  \mathbf{F}\left(\mathbf{X}^{l}\right) 
    &=& \left.\left(\mathbf{U}^{d}-\mathbf{U}^{L}\right)\right/\Delta t, \label{eq:lforce}  \\
  \mathbf{f}\left(\mathbf{x}_k\right)
    &=& \sum_{l=1}^{nl} \, c_l \, \phi_{k}^{l}(\mathbf{x}_k,\mathbf{X}^l) \, 
        \mathbf{F} \left(\mathbf{X}^{l}\right).\label{eq:spread}
\end{eqnarray}
Here, $\mathbf{U}^L$ is the interpolated velocity at the Lagrangian marker, and $\mathbf{U}^d$ is the desired velocity at the immersed boundary. $\tilde{\mathbf{u}}$ is the intermediate velocity during the time integration by only considering the nonlinear term and the viscous term in (\ref{eq:ns}). In the spreading step Eq.~(\ref{eq:spread}), the coefficient $c_l$ is calculated as
\begin{equation}\label{eq:cl}
  c_l = \frac{\Delta V^l}{\sum_{k=1}^{ne} \, \phi_{k}^l \, \Delta V_k},
\end{equation} 
where $\Delta V_k$ is the volume of the Eulerian cell $k$, and $\Delta V^l$ is the volume associated with the Lagrangian marker $l$, respectively. A common choice for the latter is $\Delta V^l = A^l h^l$ with $A^l$ being the area of the triangular element related to the Lagrangian marker $l$ and $h^l=1/3\sum_{k=1}^{ne} \phi_k^l (\Delta x_k + \Delta y_k + \Delta z_k)$~\cite{spandan2017}. Here $\Delta x$, $\Delta y$, and $\Delta z$ are the grid spaces in the three directions, respectively. For uniform Eulerian grids, $h^l$ equals to the mesh size. $c_l$ defined in Eq.~(\ref{eq:cl}) ensures the total momentum and torque conserved during the spreading step.

The transfer function $\phi^l_k$ is constructed in the same way as in Refs.~\cite{spandan2017,vanella2009}. $\phi_k^l$ can be expressed in the form of a column vector $\mathbf{\Phi}(\mathbf{X}^l)$ with length $ne$ as
\begin{equation}
   \mathbf{\Phi}^T(\mathbf{X}^l) = \mathbf{p}^T(\mathbf{X}^l)
     \mathbf{A}^{-1}(\mathbf{X}^l)\mathbf{B}(\mathbf{X}^l).
\end{equation}
Here, the superscript $T$ stands for transpose. $\mathbf{p}(\mathbf{x})$ is the column vector consisting of the linear base functions, namely $\mathbf{p}^T(\mathbf{x})=[1,~x,~y,~z]$ in the three-dimensional case. The above formula is derived by using the least-squares method with 
\begin{equation}
  \mathbf{A}(\mathbf{X}^l) = \sum_{k=1}^{ne}W(\mathbf{x}_k-\mathbf{X}^l)\mathbf{p}(\mathbf{x}_k)\mathbf{p}^T(\mathbf{x}_k),
\end{equation}
and
\begin{equation}
  \mathbf{B}(\mathbf{X}^l) = [W(\mathbf{x}_1-\mathbf{X}^l)\mathbf{p}(\mathbf{x}_1)\,,
  ~...\,,~W(\mathbf{x}_{ne}-\mathbf{X}^l)\mathbf{p}(\mathbf{x}_{ne})],
\end{equation}
in which we also use an exponential weight function 
\begin{equation}\label{eq:weight}
  W(\mathbf{x}_k-\mathbf{X}^l) = \left\{
    \begin{array}{lr}
      \exp\left[-(r_k/\alpha)^2\right], & r_k\leq 1 \\[0.2cm]
      0, & r_k> 1  
    \end{array}\right.
\end{equation}
Here $r_k=\left.\left| \mathbf{x}_k-\mathbf{X}^l \right|\,\right/H$. $\alpha$ is a constant of shape parameter. For uniform Eulerian grids one usually sets $H=1.5 \Delta x$.

\subsection{Origin of the error at the boundary in the baseline method}

Usually, the basic MLS-IBM described above could produce relatively large residual wall-normal velocity at the immersed boundary. To demonstrate this, we first notice that, after applying the IBM force, the velocity on the Eulerian grids is
\begin{equation}
  \mathbf{u}_f(\mathbf{x}_k) = \tilde{\mathbf{u}}(\mathbf{x}_k) 
  + \Delta t\,\mathbf{f}(\mathbf{x}_k),
\end{equation}
with $\mathbf{f}$ given by (\ref{eq:spread}). If we interpolate every term in the above equation back to the Lagrangian marker $\mathbf{X}^l$, we obtain
\begin{equation}
  \mathbf{U}^*(\mathbf{X}^l) = \sum_{k=1}^{ne}\phi_{k}^{l} \mathbf{u}_f^k(\mathbf{x}_k)
  = \sum_{k=1}^{ne}\phi_{k}^{l} {\bf{\tilde{u}}}(\mathbf{x}_k) 
    + \Delta t \sum_{k=1}^{ne}\phi_{k}^{l} \mathbf{f}(\mathbf{x}_k)
  = \mathbf{U}^{L}(\mathbf{X}^{l}) 
    + \Delta t \sum_{k=1}^{ne}\phi_{k}^{l} \mathbf{f}(\mathbf{x}_k),
\end{equation}
where in the last equality we use (\ref{eq:interp}). Then the actual force added to the Lagrangian marker $\mathbf{X}^l$ is
\begin{equation}
  \mathbf{F}^*(\mathbf{X}^l) = \frac{U^*(\mathbf{X}^l)-U^L(\mathbf{X}^l)}{\Delta t} 
    = \sum_{k=1}^{ne}\phi_{k}^{l} \mathbf{f}(\mathbf{x}_k)
    = \sum_{k=1}^{ne}\sum_{m=1}^{nl}c_m\phi_{k}^{l}\phi_{k}^{m}\mathbf{F}(\mathbf{X}^{m}).
\end{equation}
From this equation it is evident that the actual IBM force added onto the Lagrangian marker $\mathbf{F}^*(\mathbf{X}^l)$ is not necessarily equal to the desired value $\mathbf{F}(\mathbf{X}^l)$. This is mainly due to the fact that the IBM force at a single Eulerian grid $\mathbf{f}(\mathbf{x}_k)$ consists of the contributions from $nl$ Lagrangian markers. The fact that $\mathbf{F}^*(\mathbf{X}^l)\neq\mathbf{F}(\mathbf{X}^{l})$ will in turn cause error in the velocity boundary condition at immersed boundary.

Such phenomenon is not limited to the MLS kernel function, but also exists for the transfer functions of other types. To reduce the error, some corrections have to be made to the force $\mathbf{f}(\mathbf{x})$ so that at Lagrangian markers the actual force added is close to the desired value. Kempe~\cite{kempe2012} discussed in details about this phenomenon and proposed an iterative method to reduce the error. Another effective method to improve the accuracy is to solve a linear system for the Eulerian IBM force $\mathbf{f}$ at all the related grid points $\mathbf{x}_k$ which ensure the desired value of $\mathbf{F}$ simultaneously at all Lagrangian markers $\mathbf{X}^l$. The former requires extra computing time. While the latter is an implicit method and evolves solving a linear system with the size of the number of Lagrangian markers.

\section{The improved forcing scheme}

We now propose our new explicit and non-iterative method to reduce the error at the immersed boundary. We first made a key observation about the boundary error associated to the baseline method, and then describe the improved forcing scheme based on this key observation.

\subsection{Distribution of the values of boundary error}\label{sec:errana}

\begin{figure}[h]
  \centering
  \includegraphics[width=0.5\textwidth]{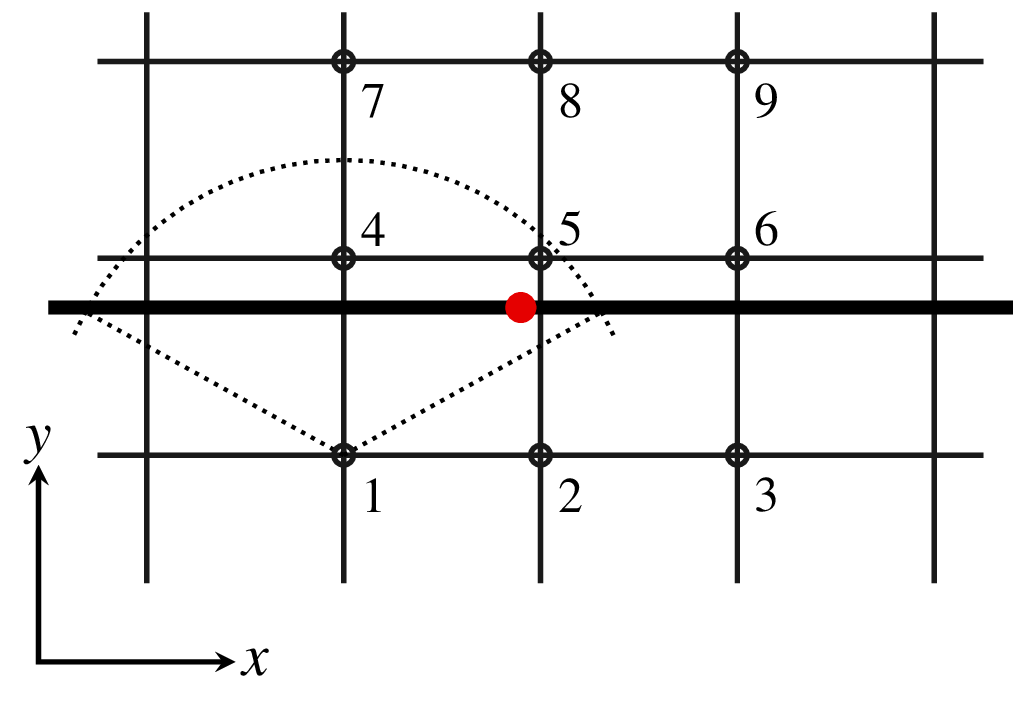}%
  \caption{The configuration of a straight immersed boundary (thick solid line) parallel to the $x$-axis and the uniform Eulerian grids (thin solid line). The red dot denotes a Lagrangian marker with coordinates $(X^0,Y^0)$. The interpolation stencil on the Eulerian grids for the Lagrangian marker $(X^0,Y^0)$ is indicated by the open circles label from 1 to 9. The dashed arc represents the region within which the Lagrangian markers contribute to grid point 1 during the spreading operation.}
  \label{fig:lerror}
\end{figure}

To analytically demonstrate the property of the error in the baseline MLS-IBM, we first look at a simple 2D situation as shown in figure~\ref{fig:lerror}. The immersed boundary is a infinite straight line parallel to the $x$-axis. Uniform Eulerian grids are used with $\Delta x=\Delta y=h$. Along this straight line we evenly distribute $n$ Lagrangian markers in each Eulerian cell. That is, each Lagrangian marker represents a line segment with the length $h/n$. To further simply the analysis, we assume the desired IBM forcing $F$ is constant along the straight line and $n$ is large enough, namely, the Lagrangian markers are dense on the straight boundary. We choose one Lagrangian marker $(X^0,Y^0)$ which is marked by the red dot, and its related Eulerian stencil grids are labelled from 1 to 9. These nine points are used for interpolating quantity at $(X^0, Y^0)$.

We now calculate the error of IBM force at the Lagrangian marker. We first determine the force at Eulerian grids. Taking point 1 for example, the IBM force at this Eulerian point is
\begin{equation}
  f_1 = \lim\limits_{n\to +\infty}\sum_{l=1}^{nl} c_l\, \phi_k^l(\mathbf{x}_1,\mathbf{X}^l)\, F 
      = \left(\lim\limits_{n\to +\infty}\sum_{l=1}^{nl} \frac{h}{hn}\, \phi_k^l(\mathbf{x}_1,\mathbf{X}^l)\right)\,F
      = \left(\lim\limits_{n\to +\infty}\sum_{l=1}^{nl} \frac{1}{n}\, \phi_k^l(\mathbf{x}_1,\mathbf{X}^l)\right)\,F
      =K_1\,F,
\label{eq:f1}
\end{equation}
where $\mathbf{X}^l$ with $l=1,...,nl$ are all the Lagrangian markers which contribute to point 1 during the spreading operation. For the current setting, they are the markers at the segment of the straight line within the circle centring at point 1 with radius $H=1.5h$, as shown in figure~\ref{fig:lerror}. Clearly, for grid points 1, 2, and 3 the coefficient $K_1$ should be same, i.e. $f_1=f_2=f_3=K_1 F$, since all three points involve the segment of the straight line with the same length, and the Lagrangian markers are evenly distributed. Actually, $K_1$ should only depends on the normal distance between the Eulerian points and the straight lines. Therefore, one has $f_4=f_5=f_6=K_4 F$, and $f_7=f_8=f_9=K_7 F$.

We then compute the actual force applied to a single Lagrangian marker after applying the IBM force at the Eulerian points. For an arbitrary Lagrangian marker indicated by the red dot in figure~\ref{fig:lerror}, the actual force $F^*$ is given by the interpolation based on the Eulerian grid points 1 to 9, saying,
\begin{equation}
  F^*=\sum_{i=1}^{9}\phi(\mathbf{x}_i,\mathbf{X}^0)f_i.
  \label{eq:Fstart}
\end{equation}
Then the ratio $F/F^*$ depends on the specific form of $\phi_k^l$ and the coordinate $Y^0$ of the considered marker $(X^0,Y^0)$. If one uses the zeroth-order constant base functions $\mathbf{p}$ and the exponential weight function $W(\mathbf{x}_k-\mathbf{X}^l)=\exp[-( r_k/\alpha)^2]$, it can be proved that 
\begin{equation}\label{eq:Fratio}
  F/F^*=1/\left(C_1K_1+C_4K_4+C_7K_7\right),
\end{equation}
with 
\begin{eqnarray*}
  C_1&=&\left[1+\exp\left(\frac{(y_1-Y^0)^2-(y_4-Y^0)^2}{\alpha^2}\right)
             +\exp\left(\frac{(y_1-Y^0)^2-(y_7-Y^0)^2}{\alpha^2}\right)\right]^{-1},\\
  C_4&=&\left[1+\exp\left(\frac{(y_4-Y^0)^2-(y_1-Y^0)^2}{\alpha^2}\right)
             +\exp\left(\frac{(y_4-Y^0)^2-(y_7-Y^0)^2}{\alpha^2}\right)\right]^{-1},\\
  C_7&=&\left[1+\exp\left(\frac{(y_7-Y^0)^2-(y_1-Y^0)^2}{\alpha^2}\right)
             +\exp\left(\frac{(y_7-Y^0)^2-(y_4-Y^0)^2}{\alpha^2}\right)\right]^{-1}.\\
\end{eqnarray*}
The details of proof are given in the Appendix. Thus, for the current configuration, the ratio $F/F^*$ is constant along the straight boundary and depends on $\alpha$ and $Y^0$, but does not on the Eulerian cell size $h$.

\begin{figure}
	\centerline{\includegraphics[width=\textwidth]{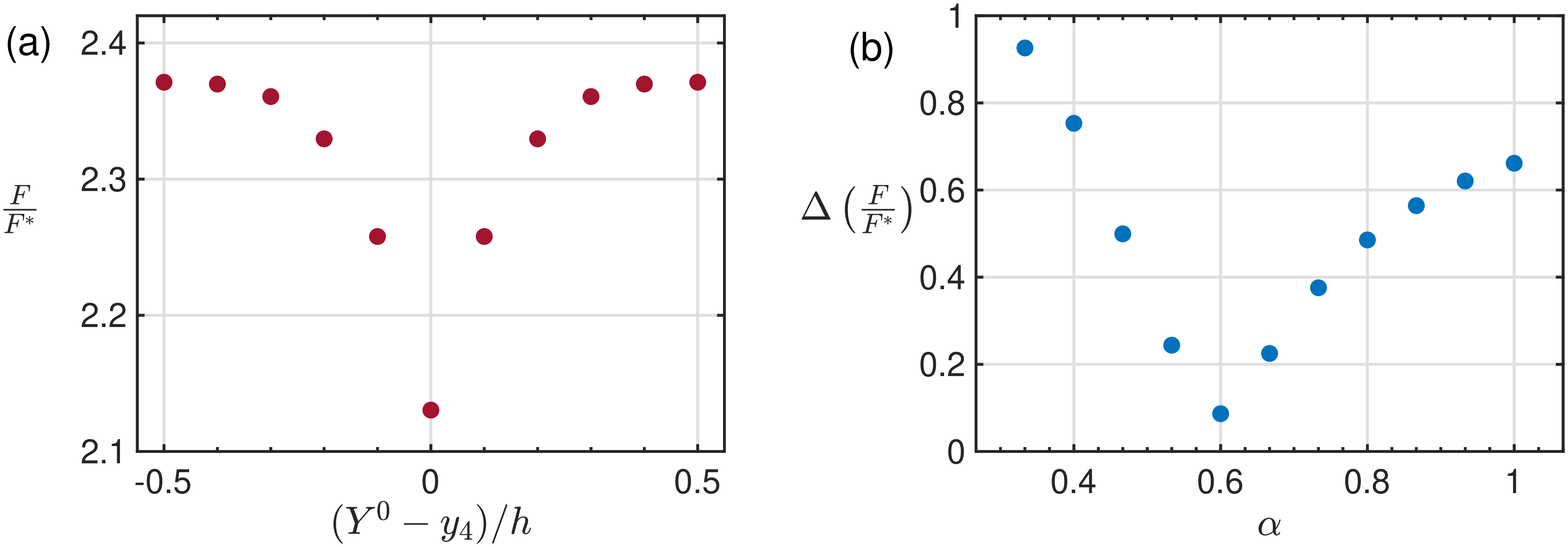}}%
    \caption{(a) The dependence of ratio $F/F^*$ on $Y^0$ for $\alpha=2/3$. (b) The maximal variation $\Delta(\frac{F}{F^*})$ for different empirical constant $\alpha$.}%
	\label{fig:Fratio}
\end{figure}

In the Appendix, we also prove that for the first-order linear base functions $\mathbf{p}^T=[1,x,y]$ and the exponential weight function $W(\mathbf{x}_k-\mathbf{X}^l)=\exp[-( r_k/\alpha)^2]$ the same observation still holds. That is, once $Y^0$ of the horizontal straight boundary is fixed, the ratio $F/F^*$ is same for all the Lagrangian markers on the boundary. The analytical expression of $F/F^*$ is rather complicated, and we test the observation by numerical experiments. The results are shown in figure~\ref{fig:Fratio}. In figure~\ref{fig:Fratio}(a) we fix $\alpha=2/3$ and move the straight boundary along the $y$-direction over a distance of $h$, saying from the midpoint of points 1 and 4 to that of points 4 and 7. The ratio $F/F^*$ is the smallest when the straight boundary coincide with point 4 or the Eulerian grid line, and the largest when it locates right at the middle of the Eulerian cell. Namely, the value is determined by the relative location of the straight boundary to the Eulerian grid lines. The difference between the largest and the smallest values, which is denoted by $\Delta(\frac{F}{F^*})$, is about $0.24$ for $\alpha=2/3$. In figure~\ref{fig:Fratio}(b) we plot $\Delta(\frac{F}{F^*})$ for different $\alpha$. Clearly, in order to minimize the variation of $F/F^*$, $\alpha$ should be set around $0.6$, instead of the commonly chosen value $2/3$.

Of course, the above analysis is for a very simplified situation. The same conclusion should be valid for the vertical straight boundary parallel to the $y$-axis with constant desired IBM force. But for the vertical straight boundary it is $X^0$ that affects the value of $F/F^*$. For general immersed boundary with complex geometry and non-uniform IBM force, one of course cannot expect a constant $F/F^*$ over the boundary. However, as will shown in the validation cases below, the value of $F/F^*$ exhibits a distribution with very narrow single peak, saying for most of the Lagrangian markers the value is very close to each others.

\subsection{The improved forcing scheme}

We have shown that for the baseline MLS-IBM, the ratio of the desired force and the actual force applied to the Lagrangian markers is very close to a constant value along the immersed boundary. Then a straightforward scheme can be proposed to correct this error in the IBM force. That is, during the spreading operation of the IBM force from a Lagrangian marker to the related Eulerian grids, saying equation~\ref{eq:spread}, we introduce a correction coefficient $Z_i$ as
\begin{equation}\label{eq:fcorrect}
  f_i\left(\mathbf{x}_k\right) = Z_i \sum_{l=1}^{nl} \, c_l \, 
        \phi_{k}^{l}(\mathbf{x}_k,\mathbf{X}^l) \, 
        F_i \left(\mathbf{X}^{l}\right),
\end{equation}
where the subscript $i$ denotes the $i$-th component of the quantity. The coefficient $Z_i$ can be fixed by the least-squares method as follow. The procedure is the same for all the components and we drop the subscript $i$ in the derivation below. The total error in the $L_2$ norm over all Lagrangian markers is
\begin{equation}\label{eq:error}
  Er_{total} = \sum_{l=1}^{NL} \left[F^*(\mathbf{X}^l)-F(\mathbf{X}^l)\right]^2.
\end{equation}
Here $NL$ is the total number of the Lagrangian markers. $F^*$ is the IBM force at marker $\mathbf{X}^l$ calculated by the interpolation of $f$ given by equation (\ref{eq:fcorrect}) as
\begin{equation}
  F^{*}(\mathbf{X}^l) = Z \sum_{k=1}^{ne}\phi_{k}^{l} f(\mathbf{x}_k)
   =Z \sum_{k=1}^{ne} \sum_{m=1}^{nk} c_m\phi_{k}^{l}\phi_{k}^{m} F(\mathbf{X}^{m}).
\end{equation}
Then the total error can be expressed in a quadratic form of $Z$, saying
\begin{equation}
  Er_{total} = a_2 Z^2 + a_1 Z + a_0,
\end{equation}
with
\[
  a_2 = \sum_{l=1}^{NL}\left[\sum_{k=1}^{ne} \sum_{m=1}^{nk} 
          c_m\phi_{k}^{l}\phi_{k}^{m} F(\mathbf{X}^{m})\right]^2, 
  \quad
  a_1 = -2\sum_{l=1}^{NL}\left[\sum_{k=1}^{ne} \sum_{m=1}^{nk} 
          c_m\phi_{k}^{l}\phi_{k}^{m} F(\mathbf{X}^{m})F(\mathbf{X}^l)\right],
  \quad
  a_0 = \sum_{l=1}^{NL} F^2(\mathbf{X}^l).
\]
Then, to minimize the total error, the correction coefficient should be set as
\begin{equation}\label{eq:zcoef}
  Z = -\frac{a_1}{2a_2} 
    = \frac{ \sum_{l=1}^{NL}\left[\sum_{k=1}^{ne}\sum_{m=1}^{nk} 
              c_m\phi_{k}^{l}\phi_{k}^{m} F(\mathbf{X}^{m})F(\mathbf{X}^l)\right]}
           { \sum_{l=1}^{NL}\left[\sum_{k=1}^{ne}\sum_{m=1}^{nk}
              c_m\phi_{k}^{l}\phi_{k}^{m} F(\mathbf{X}^{m})\right]^2}.
\end{equation}

Although the expression of the correction coefficient, namely equation~(\ref{eq:zcoef}) looks quite complex, but it only requires algebraic operations and can be readily evaluated after the desired IBM force $F(\mathbf{X}^l)$ is determined for all the Lagrangian markers during each time step in simulation. Furthermore, the coefficient only involves the transfer function $\phi$ but does not depend on the scheme of the flow solver. Therefore, the forcing correction scheme proposed here can be readily implemented to other MLS-IBM codes.

\section{Validations}\label{sec:cases}

In this section we present a series of test problems to validate the accuracy and efficiency of the proposed method. We start with the flow over a stationary sphere. Then we test our method for the moving boundary, saying oscillating sphere. Finally we show the simulation of an aorta model with relatively complex geometry. We will compare the present method with the baseline method described in section~\ref{sec:baseline}, and the iterative method used in Ref.~\cite{spandan2017}.

\subsection{Flow over a stationary sphere}

The first problem is the incompressible flow over a stationary sphere. The Reynolds number $Re=UD/\nu$ is defined with the free-stream velocity $U$, the diameter $D$, and the viscosity $\nu$. The simulation is run in a Cartesian box. In the streamwise direction a constant free-stream velocity is prescribed at the inlet, and a convective outflow condition is applied at the downstream boundary. In the two spanwise directions the periodic boundary condition is applied due to the specific configuration of our code. When the spanwise width is large enough the periodic boundary condition should have minor effect on the flow property. The Eulerian grids are always uniform with the same cell size $h$ in all three directions. On the surface the triangles associated to the Lagrangian markers always have an average edge length of about $0.7h$. Our tests show that further reducing the triangles only has very minor effect on the accuracy. The time step varies to maintain a constant CFL number of $0.2$.

\begin{figure}
	\centerline{\includegraphics[width=\textwidth]{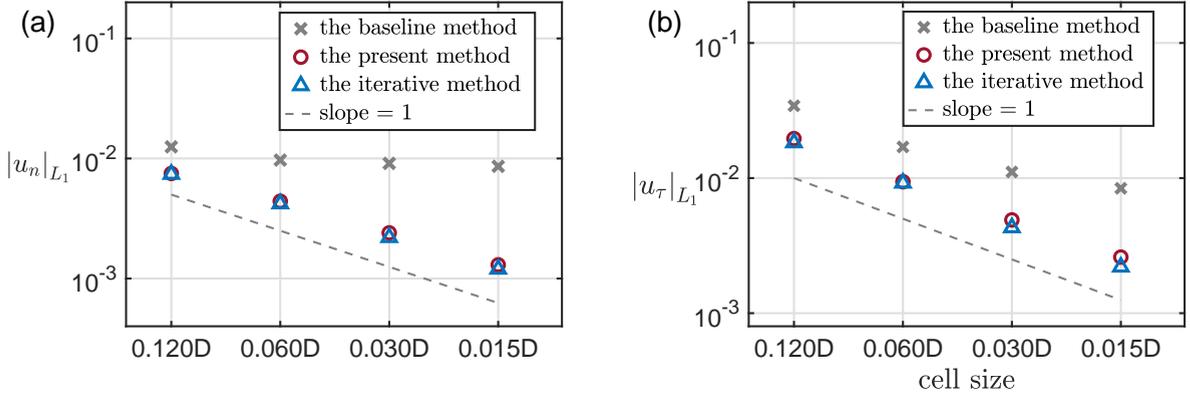}}%
	\caption{The $L_1$-norm error of the residual $u_n$ and $u_\tau$ over all Lagrangian points on a stationary sphere for different methods. The Reynolds number $Re=100$. The straight line has a slope of 1.}%
	\label{fig:erorder}
\end{figure}

We first test the order of accuracy for the boundary condition at the sphere surface for $Re=100$. The domain size is $5D \times 5D \times 5D$. We gradually reduces the Eulerian mesh size from $0.12D$ to $0.015D$ and the average size of triangles on the surface decreases accordingly. The error of the velocity boundary condition on sphere is measured by the $L_1$-norm of the residual normal velocity $u_n$ and tangential velocity $u_\tau$ over all Lagrangian markers. The results are shown in figure~\ref{fig:erorder} for three different methods, i.e. the baseline method, the present method, and the iterative method with 5 times iteration in each time step. Clearly, the baseline method has the accuracy of less than order 1. After the iterative correction the accuracy is greatly increased to order 1. The present method, explicit and without iteration, generate almost identical error level as that of the iterative method. Since only explicit operations are required, the present method cost less wall-clock time than the iterative one. For the cell size $h=0.03D$, the wall-clock time of one time step for the present method is about $25\%$ more than that of baseline method, while for the iterative method with 5 iteration the wall-clock time increases more than $100\%$ compared to the baseline method.

\begin{figure}
	\centerline{\includegraphics[width=0.6\textwidth]{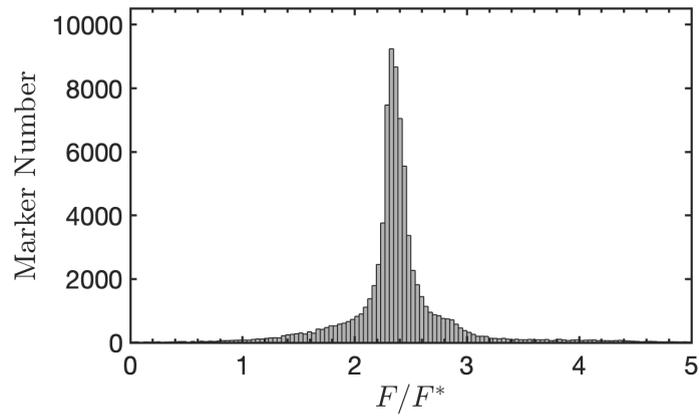}}%
	\caption{The histogram for the value of the ratio $F/F^*$ over the stationary sphere at $Re=100$. The single peak locates at $F/F^*=2.35$.}%
	\label{fig:ratiopdf}
\end{figure}

\begin{figure}
	\centerline{\includegraphics[width=0.4\textwidth]{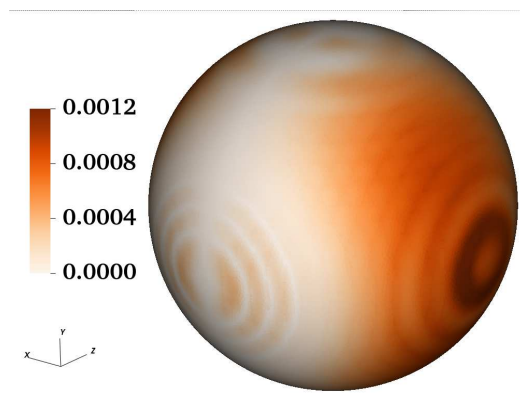}}%
	\caption{The difference between the IBM force applied in the present method and the total force applied in the iterative method. The color shows the absolute value of $F_{current}-F_{iterative}$. The free-stream flow is along the $z$-direction.}%
	\label{fig:forcediff}
\end{figure}

In the previous section we theoretically proved that for straight boundary parallel the Eulerian grid lines the ratio of the applied force $F^*$ given by the baseline method to the actually desired force $F$ is constant. For curved surface, we stated that the value of the ratio has a narrow single-peak distribution. Here we demonstrate this phenomenon for the stationary sphere by the histogram of $F/F^*$ shown in figure~\ref{fig:ratiopdf}. Indeed, a very narrow peak locates at $F/F^*=2.35$, with $79.2\%$ of total Lagrangian markers having a value in the range of $(2.0,2.7)$. We further compare the actual force applied in the present method $F_{current}$ and the total force applied in the iterative method $F_{iterative}$, saying the summation of the forces of all iterations. In figure~\ref{fig:forcediff} we plot the magnitude of the difference $F_{current}-F_{iterative}$. The maximal difference is very small, saying of the order of $10^{-3}$.

\begin{figure}
	\centerline{\includegraphics[width=\textwidth]{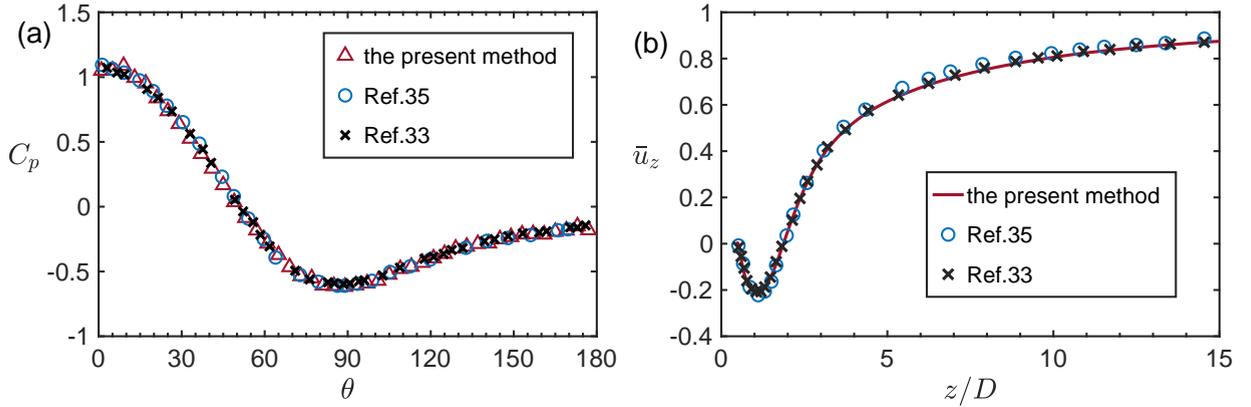}}
	\caption{(a) The pressure coefficient with azimuthal angle $\theta$ at $Re = 100$. (b) Time-averaged streamwise velocity along the z-axis at $Re = 300$.}
	\label{fig:cduwake}
\end{figure}

Having shown that the present method possesses high accuracy for the boundary condition on the immersed boundary, we now turn to the flow properties. Two Reynolds numbers are considered for $Re=100$ and $300$. For the form case the flow is steady, while for the latter one the flow is oscillating with vortex shedding. The domain size is $10D\times10D\times30D$ in the $x$, $y$, and $z$-direction, respectively. The mesh size of the uniform Eulerian grids is $h=0.03D$ for $Re=100$, and is reduced to $h=0.015D$ for $Re=300$. For both Reynolds numbers, the present method gives similar level of velocity error on the sphere surface to the iterative method. And both the present and the iterative methods predict the drag and lift coefficients, and the Strouhal number which are very close to the values given in literatures. In figure~\ref{fig:cduwake} we compare with literatures the pressure coefficient on the sphere surface for $Re=100$ in panel (a), and the time-averaged streamwise velocity in the wake along the $z$-axis for $Re=300$. Both quantities agree very well with the results reported by other groups.
\begin{table}
  \centering
  \caption{\label{tab:ssre100}Comparison of the drag coefficient $C_D$ and the $L_1$ norm of the residual velocity on the sphere among different methods and with literatures. The Reynolds number $Re = 100$.}
  \begin{tabular}{ccc}
     & $C_D$ & $\left|\mathbf{u}_{sphere}\right|_{L_1}$ \\[0.1cm]
     \hline
     the baseline method & 1.148 & 1.60e-2\\
     the iteration method & 1.117 & 3.00e-3\\
     the present method & 1.112 & 2.93e-3\\
     Wang and Zhang~\cite{wang2011} & 1.13 \\
     Kim et al.~\cite{kim2001} & 1.0875 & \\
     Fornberg~\cite{fornberg1988} & 1.085 & \\
     Johnson and Patel~\cite{johnson2000} & 1.10 &
  \end{tabular}
\end{table}
\begin{table}
  \centering
  \caption{\label{tab:ssre300}Comparison of the time averaged drag and lift coefficients, the Strouhal number, and the $L_1$ norm of the residual velocity on the sphere among different methods and with literatures. The Reynolds number $Re = 300$.}
  \begin{tabular}{ccccc}
     & $C_D$ & $C_L$ & $St$ & $\left|\mathbf{u}_{sphere}\right|_{L_1}$ \\[0.1cm]
     \hline
     the baseline method & 0.723 & 0.081 & 0.128 & 2.77e-2 \\
     the iterative method & 0.670 & 0.066 & 0.131 & 3.45e-3 \\
     the present method & 0.668 & 0.067 & 0.133 & 3.92e-3 \\
     Wang and Zhang~\cite{wang2011} & 0.680 & 0.071 & 0.135 & \\
     Kim et al.~\cite{kim2001} & 0.657 & 0.067 & 0.134 & \\
     Johnson and Patel~\cite{johnson2000} & 0.656 & 0.069 & 0.137 &\\
     Constantinescu and Squires~\cite{constantinescu2000} & 0.655 & 0.065 & 0.136 &
  \end{tabular}
\end{table}

\subsection{Flow around an oscillated sphere}

\begin{figure}
  \centering
  \includegraphics[width=0.4\textwidth]{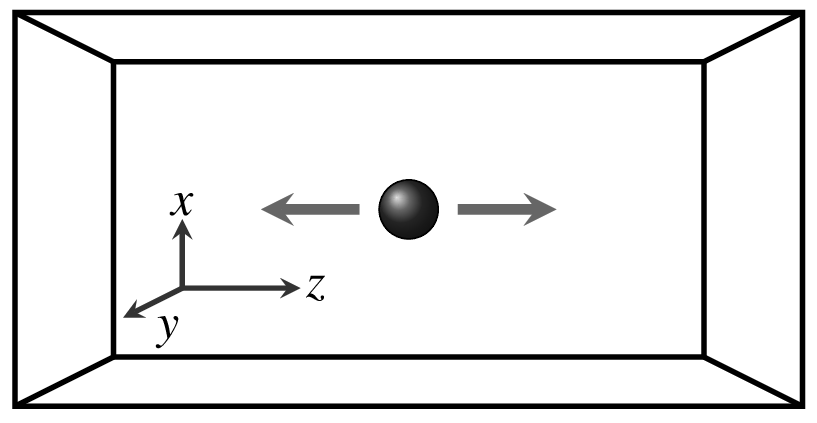}%
  \caption{The configuration for a sphere oscillating in the $z$-direction.}
  \label{fig:osconfig}
\end{figure}

The second problem is an 3D oscillating sphere of diameter $D$ in a quiescent fluid, as shown in figure~\ref{fig:osconfig}. By this test case we would like to demonstrate the performance of the present method for moving immersed boundary. The flow domain is $4D \times 4D \times 8D$. The sphere is solid and its position is given by $Z = A \sin(\omega t)$. Then the velocity of sphere is $U_z=A\omega\cos(\omega t)$. The Reynolds number is then defined as $Re=A \omega D / \nu$, i.e., by the maximal velocity and diameter of sphere. The uniform Eulerian grids have the cell size $h=0.01D$. In total 62412 triangles (or Lagrangian markers) are used to represent the sphere surface with averaged edge length of $0.0078D$. The time step varies during the simulation with the CFL number being fixed at $0.2$.

For this flow, the drag coefficient $C_D$ experiences an oscillation with the same frequency $\omega$, thus we only compare the maximal value of the drag coefficient $C_{Dmax}$ for three different amplitude $A/D$ and fixed Reynolds number $Re=100$. We also measure the accuracy of the boundary condition by calculating the $L_1$ norm of velocity error $\left|\mathbf{u}_{sphere}\right|$ at the time when the sphere reaches the maximal velocity. The results are given in table~\ref{tab:os}. Compared to the baseline method, $C_{Dmax}$ obtained by the present method is consistent with those values given in the literature. Meanwhile, the velocity error is reduced by one order of magnitude.
\begin{table}
  \centering
  \caption{\label{tab:os}Comparison of the maximal drag coefficient $C_{Dmax}$ of an oscillating sphere, and the $L_1$ norm of velocity error $\left|\mathbf{u}_{sphere}\right|$ at the time when the sphere reaches the maximal velocity. The Reynolds number is $Re = 100$.}
  \begin{tabular}{ccc}
    $A/D=0.5$ & $C_{Dmax}$ & $\left|\mathbf{u}_{sphere}\right|$ \\[0.1cm]
	\hline
    the baseline method & 3.62 & 5.22e-2\\
    the present method & 3.05 & 7.76e-3\\
	Rafi Sela et al.~\cite{sela2021} & 3.18 &\\
	Blckburn~\cite{blackburn2002} & 2.97 & \\[0.3cm]
	
    $A/D=1.0$ & $C_{Dmax}$ & $\left|\mathbf{u}_{sphere}\right|$ \\[0.1cm]
    \hline
    the baseline method & 2.45 & 8.86e-2 \\
    the present method & 2.12 & 4.13e-3 \\
    Rafi Sela et al.~\cite{sela2021} & 2.13 & \\
    Blckburn~\cite{blackburn2002} & 2.06 & \\[0.3cm]
		
    $A/D=1.5$ & $C_{Dmax}$ & $\left|\mathbf{u}_{sphere}\right|$ \\[0.1cm]
    \hline
    the baseline method & 2.05 & 8.76e-2 \\
    the present method & 1.81 & 5.20e-3\\
    Rafi Sela et al.~\cite{sela2021} & 1.82 & \\
    Blckburn~\cite{blackburn2002} & 1.77 & 
  \end{tabular}
\end{table}

\subsection{Flow inside a model aorta}

\begin{figure}
  \centering
  \includegraphics[height=5cm]{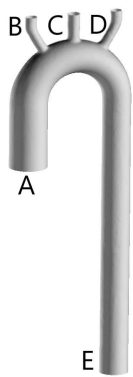}
  \caption{The geometry of the vascular model. The inlet A and outlets B,C,D have the fixed flow rate. The outlet E is left free.}
  \label{fig:aortageo}
\end{figure}

\begin{figure}
  \centering
  \includegraphics[width=\textwidth]{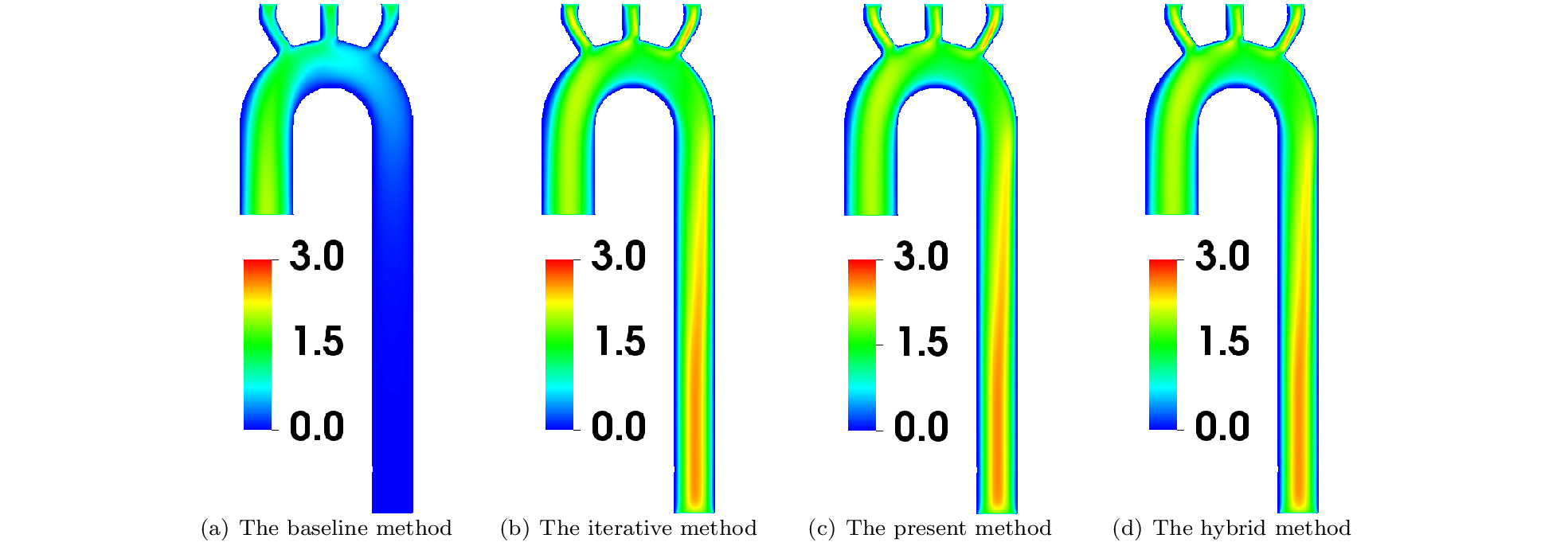}
  \caption{The contours of velocity magnitude on the vertical mid-plane for four different methods.}
  \label{fig:aortavel}
\end{figure}

For the third test case we employ a model aorta consisting of circular pipes, as shown in figure~\ref{fig:aortageo}. By this problem we aim to test the performance of the present method on the more complex geometry. The vascular model is set in the $2D \times 4D \times 10D$ rectangular region and the upper three outlets B, C, D and the lower one E are close to the boundary. The domain is periodic in the $x$ and $y$ direction and  a convective boundary condition $( \partial u_i /\partial t + c\partial u_i /\partial x = 0 )$ is used for the outflow boundary in the $z$ direction, where $c$ is the space-averaged streamwise velocity at the exit. At inlet A we prescribe a steady Poiseuille profile with maximal velocity $U_{max}$ at the central line. The Reynolds number $Re=U_{max}D_A/\nu$ is set at $300$. Here $D_A$ is the diameter of inlet A. The diameters for the three top outlets B, C, and D are $D_A/3$, and that for the bottom outlet E is $4D_A/5$, respectively. The flow rates at B, C, and D are all kept constant and equal to $5\%$ of the flow rate at the inlet A. The outlet E is let free without any control of flow rate. Since the outlet E is very close to the lower boundary set as a convective boundary condition, the flow slows down as approaching the outlet E, as can be seen from figure~\ref{fig:aortavel}.

The aorta is placed inside a Cartesian box with uniform Eulerian grids. The flow field grid is $128 \times 256 \times 640$ and the Euler grid size is $0.015D$. There are a total of 463051 Lagrangian markers, with an average edge length of $0.01D$. All simulations are performed by using one whole AMD EPYC processor (model 7452 with 64 cores and 256 GB RAM).

For this flow we compare the baseline method, the iterative method, and the present method. For the iterative method we use 10 times of iteration for each time step. We also test the hybrid method, in which the new proposed force correction procedure is conducted at the end of each iteration in the iterative method. Such implement is straightforward due to the explicit nature of our method. In the hybrid method, only two times of iterations are carried out for each time step. In figure~\ref{fig:aortavel} we plot the contours of velocity magnitude on the vertical mid plane for the four methods. Clearly, the baseline method suffers from severe mass loss through the aorta boundary. There is nearly no vertical flow in the right main pipe. For the other three methods, strong downward flows develop in the vertical pipe on the right. 

\begin{figure}
  \centering
  \includegraphics[width=\textwidth]{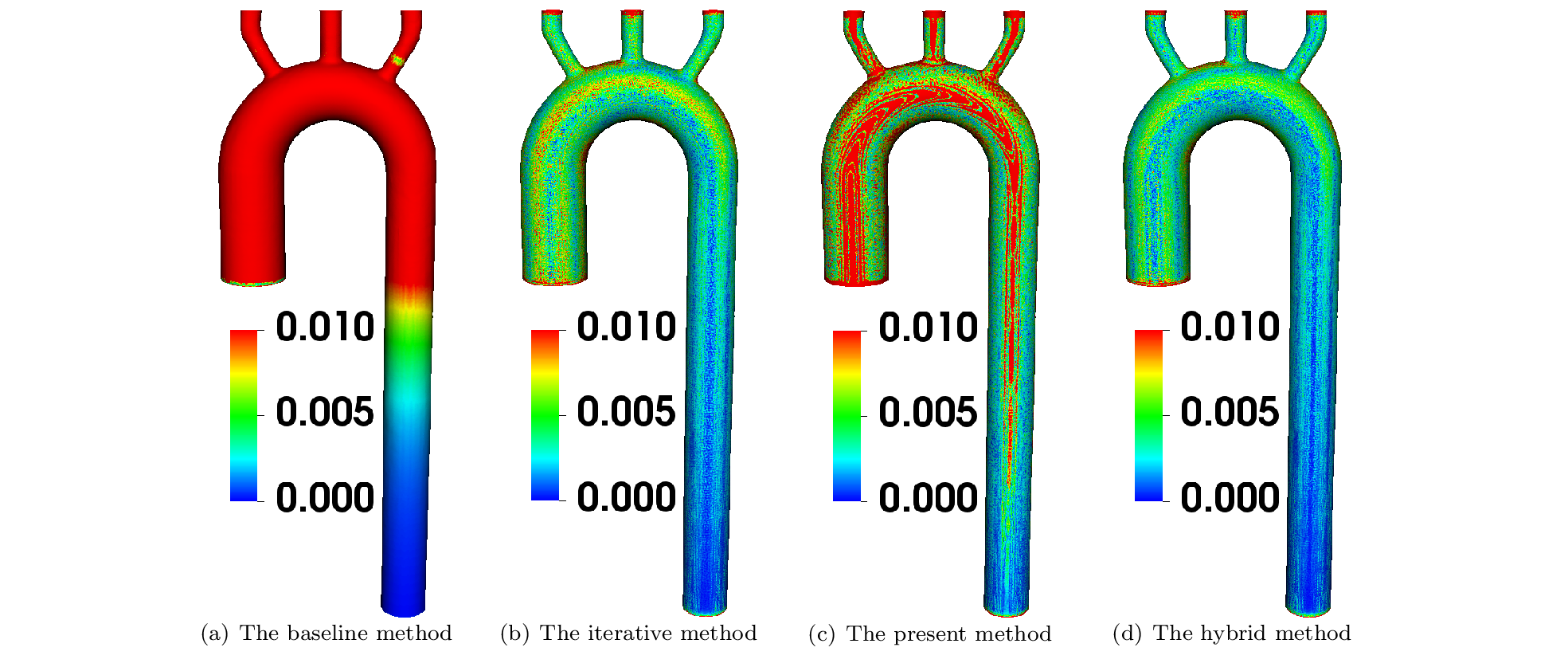}
  \caption{The error of velocity magnitude on the aorta boundary for four different methods.}
  \label{fig:aortaerr}
\end{figure}

The $L_1$ norm of the wall velocity is calculated and summarised in table~\ref{tab:aorta}. We did not show the error for the baseline method since the flow in the right pipe does not develop. The present method along produces an error about twice the value of the iterative method with 10 times iteration. The hybrid method, which consists of 2 times iterations at each time step, already has an error smaller  than the pure iterative method. Of course, less iterations correspond to less wall-clock time, as compared in table~\ref{tab:aorta}. We further display the distributions of the wall-velocity magnitude for the four methods in figure~\ref{fig:aortaerr}.
\begin{table}
  \centering
  \caption{\label{tab:aorta}The averaged residual Velocity on the surface and the averaged wall-clock time for one step of integration for different methods.}
  \begin{tabular}{ccc}
    & The averaged residual Velocity$\overline{|u|}$& the averaged wall-clock time(s)  \\
    \hline
    the baseline method & $-$ & 2.08s \\
	the iterative method & 3.23e-3 & 14.63s \\
    the present method & 7.73e-3 & 3.21s \\
	the hybrid method & 3.05e-3 & 4.85s \\
	\end{tabular}
\end{table}

\section{Conclusions}\label{sec:conclusion}

In summary, we develop an explicit and non-iterative technique to construct a MLS-IBM method with low error for the velocity condition on the immersed boundary. The technique bases on a key observation about the value distribution of the ratio of the actually applied IBM force on the Lagrangian markers to the desired value. For straight boundary parallel to the Eulerian grid lines, the ratio is constant along the boundary and depends on the relative location of the boundary to the grid lines. For immersed boundary with complex geometry, numerical results reveal that the value of the ratio exhibits a narrow peak distribution, saying for most Lagrangian markers the ratio is very close to certain value. Based on this observation, a force correction procedure is proposed. A single coefficient is introduced to correct the IBM force, and its value can be readily determined by a least-squares method. This correction is explicit and non-iterative.

The new method is tested for the flow over a stationary sphere, the flow around an oscillating sphere with moving boundary, and the flow inside a model aorta with complex boundary geometry. The performance is compared with the baseline MLS-IBM method and the iterative method. Simulation results suggest that the present method can achieve the low level of boundary error similar to the iterative method, but requires less computing time. The new technique can be readily combined with the iterative method so that the boundary error can be further reduced with less iteration times compared to the purely iterative method.

The force correction method proposed here can be easily implemented into other MLS-IBM code since it only involves algebraic calculation. Although the current method is developed specific for the transfer function used in MLS-IBM, but for other transfer function, the same strategy may still be applicable. It is therefore of interests to test the current method in IBM method other than the MLS type. 

\section*{Acknowledgements}
The support from the Major Research Plan of National Natural and Science Foundation of China for Turbulent Structures under Grants 91852107 and 91752202.

\section*{Appendix}

In section~\ref{sec:errana} we stated that for the simple 2D straight boundary as shown in figure~\ref{fig:lerror}, the ratio $F/F^*$ is constant and independent of the Eulerian cell size $h$ for the zeroth-order and first-order basis functions $\mathbf{p}$. In this appendix we provide the details of the proof. 

For the zeroth-order constant basis function, MLS interpolation degenerates into shepard-interpolation, and the transfer function $\phi_{k}^{l}$ is:
\begin{equation}
\phi_{k}^{l}(\mathbf{x}_k, \mathbf{X}^l)
=\frac{W(\mathbf{x}_k-\mathbf{X}^l)}{\sum^{ne}_{m=1}W(\mathbf{x}_m-\mathbf{X}^l) }
=\frac{e^{-( r_k/\alpha)^2}}{\sum_{m=1}^{9}e^{-( r_m/\alpha)^2} }.
\end{equation}
The summation of the MLS transfer functions on points $\mathbf{x}_{i=1,2,3}$ is
\begin{equation}
\sum_{i=1}^{3}\phi(\mathbf{x}_i,\mathbf{X}^0)
=\frac{\sum_{k=1}^{3}e^{-(r_k/\alpha)^2}}
{\sum^{9}_{m=1}e^{-( r_m/\alpha)^2} }.
\label{eq:phi1}
\end{equation}
This relation can be simplified by noticing that
\begin{eqnarray*}
\frac{\sum_{k=1}^{3}e^{-( r_k/\alpha)^2}}
{\sum^{6}_{m=4}e^{-( r_m/\alpha)^2} }
&=&\frac{\sum_{k=1}^{3}e^{\frac{-[(x_k-X^0)^2+(y_k-Y^0)^2]}{\alpha ^2}}}
{\sum^{6}_{m=4}e^{\frac{-[(x_m-X^0)^2+(y_m-Y^0)^2]}{\alpha ^2}} }
=\frac{e^{\frac{-(y_1-Y^0)^2}{\alpha^2}}}
{e^{\frac{-(y_4-Y^0)^2}{\alpha^2}}}
=e^{\frac{-(y_1-Y^0)^2+(y_4-Y^0)^2}{\alpha^2}},\\
\frac{\sum_{k=1}^{3}e^{-( r_k/\alpha)^2}}
{\sum^{9}_{m=7}e^{-( r_m/\alpha)^2} }
&=&\frac{\sum_{k=1}^{3}e^{\frac{-[(x_k-X^0)^2+(y_k-Y^0)^2]}{\alpha ^2}}}
{\sum^{9}_{m=7}e^{\frac{-[(x_m-X^0)^2+(y_m-Y^0)^2]}{\alpha ^2}} }
=\frac{e^{\frac{-(y_1-Y^0)^2}{\alpha^2}}}
{e^{\frac{-(y_7-Y^0)^2}{\alpha^2}}}
=e^{\frac{-(y_1-Y^0)^2+(y_7-Y^0)^2}{\alpha^2}}.
\end{eqnarray*}
Substituting the above two equations back to (\ref{eq:phi1}), one obtains
\begin{equation}
\sum_{i=1}^{3}\phi(\mathbf{x}_i,\mathbf{X}^0)
=C_1
\quad with \quad C_1=\left[1+\exp\left(\frac{(y_1-Y^0)^2-(y_4-Y^0)^2}{\alpha^2}\right)
+\exp\left(\frac{(y_1-Y^0)^2-(y_7-Y^0)^2}{\alpha^2}\right)\right]^{-1}
\label{eq:c1}
\end{equation}
Similar relations can be derived for the MLS transfer functions on points $\mathbf{x}_{i=4,5,6}$ and $\mathbf{x}_{i=7,8,9}$, respectively, as
\begin{eqnarray}
\sum_{i=4}^{6}\phi(\mathbf{x}_i,\mathbf{X}^0) = C_4
\quad with \quad C_4=\left[1+\exp\left(\frac{(y_4-Y^0)^2-(y_1-Y^0)^2}{\alpha^2}\right)
+\exp\left(\frac{(y_4-Y^0)^2-(y_7-Y^0)^2}{\alpha^2}\right)\right]^{-1} \label{eq:c2}\\
\sum_{i=7}^{9}\phi(\mathbf{x}_i,\mathbf{X}^0) =C_7
\quad with \quad C_7=\left[1+\exp\left(\frac{(y_7-Y^0)^2-(y_1-Y^0)^2}{\alpha^2}\right)
+\exp\left(\frac{(y_7-Y^0)^2-(y_4-Y^0)^2}{\alpha^2}\right)\right]^{-1}
\label{eq:c3}
\end{eqnarray}
Clearly, $C_{i=1,4,7}$ are independent of $X^0$. Then (\ref{eq:Fstart}) and (\ref{eq:Fratio}) can be readily obtained, and $F/F^*$ is independent of $X^0$. Moreover, since $\phi_{k}^{l}(\mathbf{x}_k, \mathbf{X}^l) = \exp[-( r_k/\alpha)^2]/\sum_{m=1}^{9}\exp[-( r_m/\alpha)^2]$ is function of the normalized local coordinates as $r_k=\left.\left| \mathbf{x}_k-\mathbf{X}^l \right|\,\right/H$ with $H=1.5 h$, $F/F^*$ is also independent of $h$, i.e. only depending on the relative location with the grid line.

For the first-order linear basis functions, the same two conclusions still hold. The MLS transfer function $\phi_k^l$ can be expressed in the form of a column vector $\mathbf{\Phi}(\mathbf{X}^l)$ with length $ne$ as
\begin{equation}
  \mathbf{ \Phi }^{T}(\mathbf{X}^l) 
   = \mathbf{p}^T(\mathbf{X}^l)\mathbf{A}^{-1}(\mathbf{X}^l) \mathbf{B}(\mathbf{X}^l)
   = [ \phi(\mathbf{x}_1,\mathbf{X}^0)~...~\phi(\mathbf{x}_9,\mathbf{X}^0) ]_{1\times9}.
\label{eq:Phi}
\end{equation}
Here $\mathbf{p}^T(\mathbf{x})$ represents the one-order linear basis function $\mathbf{p}^T(\mathbf{x})=[1, x, y]$. And the two matrices are
\[
  \mathbf{A}(\mathbf{X}^l)
    =\sum_{k=1}^{ne}W(\mathbf{x}_k-\mathbf{X}^l)\mathbf{p}(\mathbf{x}_k)\mathbf{p}^T(\mathbf{x}_k)
    ={ \left[ \begin{array}{ccc}
	   \sum_{i=1}^{9}W_i & \sum_{i=1}^{9}x_iW_i & \sum_{i=1}^{9}y_iW_i\\
	   \sum_{i=1}^{9}x_iW_i & \sum_{i=1}^{9}{x_i}^2W_i & \sum_{i=1}^{9}x_iy_iW_i\\
	   \sum_{i=1}^{9}y_iW_i & \sum_{i=1}^{9}x_iy_iW_i & \sum_{i=1}^{9}{y_i}^2W_i\\
	   \end{array} 
	 \right] }_{(3\times 3)},
\]
and
\[
  \mathbf{B}(\mathbf{X}^l)
    =[W(\mathbf{x}_1-\mathbf{X}^l)\mathbf{p}(\mathbf{x}_1)
        \,,~...\,,W(\mathbf{x}_{ne}-\mathbf{X}^l)\mathbf{p}(\mathbf{x}_{ne})]
    ={ \left[ \begin{array}{cccc}
	     W_1 & W_2 & ...&W_9\\
	     x_1W_1 & x_2W_2 & ...&x_9W_9\\
	     y_1W_1 & y_2W_2 & ...&y_9W_9\\
	    \end{array} 
	   \right]}_{(3\times 9)}.
\]
We further denote the inverse of $\mathbf{A}$ as
\begin{equation}
\mathbf{A}^{-1}(\mathbf{X}^l)={
	\left[ \begin{array}{ccc}
	a_{11} & a_{12} & a_{13}\\
	a_{21} & a_{22} & a_{23}\\
	a_{31} & a_{32} & a_{33}\\
	\end{array} 
	\right ]}_{(3\times 3)}.
\end{equation}

Like our proof for the zeroth-order basic function, in order to have $F/F^*$ independent of $X^0$, $\sum_{i=1}^{3}\phi(\mathbf{x}_i,\mathbf{X}^0)$, $\sum_{i=4}^{6}\phi(\mathbf{x}_i,\mathbf{X}^0)$, and $\sum_{i=7}^{9}\phi(\mathbf{x}_i,\mathbf{X}^0)$ should all be independent of $X^0$. We give the procedure for $\sum_{i=1}^{3}\phi(\mathbf{x}_i,\mathbf{X}^0)$, and the other two can be treated similarly. From the definitions, $\phi(\mathbf{x}_i,\mathbf{X}^0)$ for $i=1, 2, 3$ read
\[
\phi(\mathbf{x}_1,\mathbf{X}^0)
=(a_{11}+a_{12}X^0+a_{13}Y^0)W_1
+(a_{21}+a_{22}X^0+a_{23}Y^0)x_1W_1
+(a_{31}+a_{32}X^0+a_{33}Y^0)y_1W_1,
\]
\[
\phi(\mathbf{x}_2,\mathbf{X}^0)
=(a_{11}+a_{12}X^0+a_{13}Y^0)W_2
+(a_{21}+a_{22}X^0+a_{23}Y^0)x_2W_2
+(a_{31}+a_{32}X^0+a_{33}Y^0)y_2W_2,
\]
\[
\phi(\mathbf{x}_3,\mathbf{X}^0)
=(a_{11}+a_{12}X^0+a_{13}Y^0)W_3
+(a_{21}+a_{22}X^0+a_{23}Y^0)x_3W_3
+(a_{31}+a_{32}X^0+a_{33}Y^0)y_3W_3.
\]
The summation of the above three equations gives, and by using the symmetry of $\mathbf{A}^{-1}$,
\begin{eqnarray}
  \sum_{i=1}^{3}\phi(\mathbf{x}_i,\mathbf{X}^0)
   &=&[a_{11}\sum_{i=1}^{3}W_i+a_{12}\sum_{i=1}^{3}x_iW_i+a_{13}\sum_{i=1}^{3}y_iW_i]\nonumber\\
   &+&X^0[a_{21}\sum_{i=1}^{3}W_i+a_{22}\sum_{i=1}^{3}x_iW_i+a_{23}\sum_{i=1}^{3}y_iW_i]\nonumber\\
   &+&Y^0[a_{31}\sum_{i=1}^{3}W_i+a_{32}\sum_{i=1}^{3}x_iW_i+a_{33}\sum_{i=1}^{3}y_iW_i].
   \label{eq:phic}
\end{eqnarray}

We first prove that $a_{23}=a_{32}=0$. The inverse of $\mathbf{A}$ can be calculated as $\mathbf{A}^{-1}=\mathbf{A}^*/|\mathbf{A}|$ with $\mathbf{A}^*$ being the adjoint matrix and $|\mathbf{A}|$ the determinant, respectively. Then the explicit expression of $a_{23}$ reads
\begin{eqnarray}
    a_{23} &=& -\frac{1}{|\mathbf{A}|}\left(\sum_{i=1}^{9}W_i \, \sum_{i=1}^{9}x_iy_iW_i)
           - \sum_{i=1}^{9}x_iW_i \, \sum_{i=1}^{9}y_iW_i\right) \nonumber \\
    &=& -\frac{h^2}{|\mathbf{A}|}\left[~(W_1-W_3-W_7+W_9)\right.
                        \left(\sum_{i=1}^{9}W_i\right) \nonumber \\
     &&\quad\quad\quad \left.-(W_3+W_6+W_9-W_1-W_4-W_7)(W_7+W_8+W_9-W_1-W_2-W_3)~\right].
\label{eq:a23}
\end{eqnarray}
By noticing that 
\[
  \frac{W_7}{W_1}=\frac{W_8}{W_2}=\frac{W_9}{W_3}
      =\frac{\sum_{i=7}^{9}W_i}{\sum_{i=1}^{3}W_i}=R_{7,1},
  \quad\quad
  \frac{W_4}{W_1}=\frac{W_5}{W_2}=\frac{W_6}{W_3}
      =\frac{\sum_{i=4}^{6}W_i}{\sum_{i=1}^{3}W_i}=R_{4,1},
\]
one has $\sum_{i=1}^{9}W_i=\sum_{i=1}^{3}W_i(1+R_{4,1}+R_{7,1})$ and equation~(\ref{eq:a23}) gives 
\begin{eqnarray}
  a_{23} &=& -\frac{h^2}{|\mathbf{A}|}
         \left[(W_1-W_3)(1-R_{7,1})(\sum_{i=1}^{3}W_i)(1+R_{4,1}+R_{7,1})\right.\nonumber\\
         && \quad \quad \quad
         \left. -(W_1-W_3)(1+R_{4,1}+R_{7,1})(1-R_{7,1})(\sum_{i=1}^{3}W_i)\right] = 0.
\end{eqnarray}

Next, by the relation $\mathbf{A}^{-1}\mathbf{A}=\mathbf{I}$, the equalities for element $(1,1)$, $(2,1)$, and $(3,3)$ give
\begin{equation}
a_{11}\sum_{i=1}^{9}W_i+a_{12}\sum_{i=1}^{9}x_iW_i+a_{13}\sum_{i=1}^{9}y_iW_i=1.
\label{eq:a11}
\end{equation}
\begin{equation}
a_{21}\sum_{i=1}^{9}W_i+a_{22}\sum_{i=1}^{9}x_iW_i+a_{23}\sum_{i=1}^{9}y_iW_i=0.
\label{eq:a21}
\end{equation}
\begin{equation}
a_{31}y_i\sum_{i=1}^{9}W_i+a_{32}\sum_{i=1}^{9}x_iy_iW_i+a_{33}\sum_{i=1}^{9}{y_i}^2W_i=1
\label{eq:a33}
\end{equation}
Since $a_{32}=0$, equation (\ref{eq:a33}) can be cast into
\[
  a_{31}\sum_{i=1}^{9}W_i+a_{32}\sum_{i=1}^{9}x_iW_i+a_{33}\sum_{i=1}^{9}y_iW_i
  =a_{31}(1+R_{4,1}+R_{7,1})\sum_{i=1}^{3}W_i+a_{33}(y_1+y_4R_{4,1}+y_7R_{7,1})\sum_{i=1}^{3}W_i
  =0,
\]
which means $a_{31}=a_{13}=Ca_{33}$ with $C=-(y_1+y_4R_{4,1}+y_7R_{7,1})/(1+R_{4,1}+R_{7,1})$. While (\ref{eq:a11}) can be rewritten as
\[
  (1+R_{4,1}+R_{7,1})a_{11}\sum_{i=1}^{3}W_i
  +(1+R_{4,1}+R_{7,1})a_{12}\sum_{i=1}^{3}x_iW_i
  +(y_1+y_4R_{4,1}+y_7K_{7,1})a_{13}\sum_{i=1}^{3}W_i=1,
\]
or
\[
  a_{11}\sum_{i=1}^{3}W_i+a_{12}\sum_{i=1}^{3}x_iW_i
   =\frac{1}{1+R_{4,1}+R_{7,1}}+Ca_{13}\sum_{i=1}^{3}W_i.
\]
Then the first term on the right-hand-side of equation~(\ref{eq:phic}) equals to
\begin{equation}\label{eq:phic1}
    a_{11}\sum_{i=1}^{3}W_i+a_{12}\sum_{i=1}^{3}x_iW_i+a_{13}\sum_{i=1}^{3}y_iW_i
    = \frac{1}{1+R_{4,1}+R_{7,1}}+
    C(C+y_1)a_{33}\sum_{i=1}^{3}W_i
\end{equation}
Equation (\ref{eq:a21}) can be simplified as
\begin{eqnarray*}
    && a_{21}\sum_{i=1}^{9}W_i+a_{22}\sum_{i=1}^{9}x_iW_i+a_{23}\sum_{i=1}^{9}y_iW_i
    =a_{21}\sum_{i=1}^{9}W_i+a_{22}\sum_{i=1}^{9}x_iW_i \\
    && \quad\quad = (1+R_{4,1}+R_{7,1})
         \left[a_{21}\sum_{i=1}^{3}W_i+a_{22}\sum_{i=1}^{3}x_iW_i\right] = 0,
\end{eqnarray*}
which suggests that 
\[
  a_{21}\sum_{i=1}^{3}W_i+a_{22}\sum_{i=1}^{3}x_iW_i=0.
\]
Together with $a_{23}=0$, the second term on the right-hand-side of equation~(\ref{eq:phi1}) vanishes. For the third term we have, considering $a_{32}=0$,
\begin{equation}\label{eq:phic3}
Y^0\left[a_{31}\sum_{i=1}^{3}W_i+a_{32}\sum_{i=1}^{3}x_iW_i+a_{33}\sum_{i=1}^{3}y_iW_i\right]=
Y^0\left(y_1+C\right)a_{33}\sum_{i=1}^{3}W_i.
\end{equation}
Then combining equations~(\ref{eq:phic1}) and (\ref{eq:phic3}), and the fact that the second term varnishes, one finally obtains
\begin{equation}
   \sum_{i=1}^{3}\phi(\mathbf{x}_i,\mathbf{X}^0)
      = \frac{1}{1+R_{4,1}+R_{7,1}}+
      (C+Y^0)(y_1+C)a_{33}\sum_{i=1}^{3}W_i.
\end{equation}
In the above equation, the first term on the right-hand side and the prefactor $(C+Y^0)(y_1+C)$ are both independent of $X^0$. Moreover, with $a_{32}=0$, equation~(\ref{eq:a33}) gives
\begin{eqnarray*}
   && a_{31}y_i\sum_{i=1}^{9}W_i+a_{32}\sum_{i=1}^{9}x_iy_iW_i+a_{33}\sum_{i=1}^{9}{y_i}^2W_i \\
   && \quad = a_{31}[y_1+y_4R_{4,1}+y_7R_{7,1}]\sum_{i=1}^{3}W_i
             +a_{33}[{y_1}^2+{y_4}^2R_{4,1}+{y_7}^2R_{7,1}]\sum_{i=1}^{3}W_i \\
   && \quad =\left[{y_1}^2+{y_4}^2R_{4,1}+{y_7}^2R_{7,1}
                   +C(y_1+y_4R_{4,1}+y_7R_{7,1})\right]
             a_{33}\sum_{i=1}^{3}W_i=1.
\end{eqnarray*}
Again, the coefficient in the square bracket does not contain $X^0$, and the above relation implies that $a_{33}\sum_{i=1}^{3}W_i$ is independent of $X^0$. We now complete the proof that $\sum_{i=1}^{3}\phi(\mathbf{x}_i,\mathbf{X}^0)=C_1$ is independent of $X^0$. By the same way one can prove that $\sum_{i=4}^{6}\phi(\mathbf{x}_i,\mathbf{X}^0)=C_2$ and $\sum_{i=7}^{9}\phi(\mathbf{x}_i,\mathbf{X}^0)=C_3$ are both independent of $X^0$. Therefore, for the first-order basic functions, the ratio $F/F^*$ does not change with $X^0$.

Finally, we prove that the ratio $F/F^*$ is also independent of $h$ if the relative location of the straight boundary and the Eulerian grid line does not change, namely, $r_k=\left.\left| \mathbf{x}_k-\mathbf{X}^l \right|\,\right/H$ keeps the same as $h$ varies. When $h$ changes, either increases or decreases, the global coordinates for the Eulerian grid points 1 to 9 in figure~\ref{fig:lerror} changes, which corresponds to a combination of stretch and translation transformation for the basis function $\mathbf{p}$. We denote this transformation as
\[
  \mathbf{\tilde{p}}(\mathbf{X}^l) = \mathbf{C}\mathbf{p}(\mathbf{X}^l),
\]
with $\mathbf{C}$ being the invertible transformation matrix. Here the tilde denotes the quantities after the transformation. Then by definition one has
\begin{eqnarray*}
  \mathbf{\tilde{A}}(\mathbf{X}^l) 
    &=& \sum_{k=1}^{ne}W(\mathbf{x}_k-\mathbf{X}^l)
           \mathbf{\tilde{p}}(\mathbf{x}_k)\mathbf{\tilde{p}}^T(\mathbf{x}_k) 
       = \sum_{k=1}^{ne}W(\mathbf{x}_k-\mathbf{X}^l)\,\mathbf{C}
           \,\mathbf{p}(\mathbf{x}_k)\mathbf{p}^T(\mathbf{x}_k)\mathbf{C}^T, \\
  \mathbf{\tilde{B}}(\mathbf{X}^l)
    &=& [W(\mathbf{x}_1-\mathbf{X}^l)\mathbf{\tilde{p}}(\mathbf{x}_1)~...~
         W(\mathbf{x}_{ne}-\mathbf{X}^l)\mathbf{\tilde{p}}(\mathbf{x}_{ne})]
       = [W(\mathbf{x}_1-\mathbf{X}^l)\mathbf{C}\mathbf{p}(\mathbf{x}_1)~...~
          W(\mathbf{x}_{ne}-\mathbf{X}^l)\mathbf{C}\mathbf{p}(\mathbf{x}_{ne})] \\
    &=& \mathbf{C}[W(\mathbf{x}_1-\mathbf{X}^l)\mathbf{p}(\mathbf{x}_1)~...~
                   W(\mathbf{x}_{ne}-\mathbf{X}^l)\mathbf{p}(\mathbf{x}_{ne})] \\
    &=& \mathbf{C}\mathbf{B}(\mathbf{X}^l), \\
  \mathbf{ \tilde{\Phi} }^{T}(\mathbf{X}^l)
    &=& \mathbf{ \tilde{p}}^T(\mathbf{X}^l)\mathbf{\tilde{A}}^{-1}(\mathbf{X}^l)
           \mathbf{ \tilde{B}}(\mathbf{X}^l)\\
    &=& \mathbf{p}^T(\mathbf{X}^l)\mathbf{C}^T \mathbf{
           \tilde{A}}^{-1}(\mathbf{X}^l)\mathbf{C}\mathbf{B}(\mathbf{X}^l) \\
    &=& \mathbf{p}^T(\mathbf{X}^l) 
          \left[(\mathbf{C})^{-1}\mathbf{\tilde{A}}(\mathbf{X}^l)(\mathbf{C}^T)^{-1}\right]^{-1}
        \mathbf{B}(\mathbf{X}^l) \\
    &=& \mathbf{p}^T(\mathbf{X}^l) 
          \left[(\mathbf{C})^{-1}\sum_{k=1}^{ne}W(\mathbf{x}_k-\mathbf{X}^l)\mathbf{C}\,
            \mathbf{p}(\mathbf{x}_k)\mathbf{p}^T(\mathbf{x}_k)\mathbf{C}^T
            (\mathbf{C}^T)^{-1}\right]^{-1}\mathbf{B}(\mathbf{X}^l) \\
    &=& \mathbf{p}^T(\mathbf{X}^l) 
          \left[\sum_{k=1}^{ne}W(\mathbf{x}_k-\mathbf{X}^l)\mathbf{p}(\mathbf{x}_k)
           \mathbf{p}^T(\mathbf{x}_k)\right]^{-1}\mathbf{B}(\mathbf{X}^l) \\
    &=& \mathbf{ \Phi }^{T}(\mathbf{X}^l).
\end{eqnarray*}
Therefore, the transfer function $\mathbf{\Phi}(\mathbf{X}^l)$ is invariant under the transformation, and the ratio $F/F^*$ is independent of $h$.

\bibliographystyle{model1-num-names}

\end{document}